\documentstyle[12pt, amssymb, amscd]{article}

\def\ca#1{#1}
\def\nobd{}

\def\itbox#1{\ifvmode\indent\fi\makebox[1em][c]{\rm(#1)}\hskip .5em}
\def\itboxx#1{\ifvmode\indent\fi\makebox[1em][r]{\rm(#1)}\hskip .5em}
\def\mboenvskip{\ifvmode\removelastskip\fi\vskip\baselineskip}
\def\qed{\ifmmode\square\else{\noindent\unskip\nobreak\hfil
\penalty50\hskip1em\null\nobreak\hfil$\square$
\parfillskip=0pt\finalhyphendemerits=0\endgraf}\fi\mboenvskip}

\def\PfEng{\mboenvskip\noindent{\it Proof}}

\newcommand\gal{G}   %% Galois group Gal(k^s/k)
\newcommand\cok{{\cal O}_K} %% ring of integers
\newcommand\tor{T}   %% torus over \cok
\newcommand\fltor{{\cal T}}   %% formal torus over \cok
\newcommand\fla{{\cal A}}   %% formal completion of A.
\newcommand\flb{{\cal B}}   %% formal good reduction part of A.

\newcommand{\bbz}{{\Bbb Z}}
\newcommand{\ovl}{\overline}
\newcommand\Pic{\mbox{\rm Pic}}
\newcommand\chara{\mbox{\rm char}}
\newcommand\Gal{\mbox{\rm Gal}}
\newcommand\Hom{\mbox{\rm Hom}}
\newcommand\Spec{\mbox{\rm Spec\kern 2pt}}

\newcommand\Ker{\mbox{\rm Ker}}
\newcommand\Img{\mbox{\rm Im}}
\newcommand\divides{\kern 2pt|\kern 2pt}

\newtheorem{thm}{Theorem}[section]
\newtheorem{lem}[thm]{Lemma}
\newtheorem{cor}[thm]{Corollary} 
\newtheorem{prop}[thm]{Proposition} 
\newtheorem{rem}[thm]{Remark}
\newtheorem{exmp}[thm]{Example}

\begin{document}

\title{Rational points of the group of components 
of a N\'eron model}
\author{Siegfried Bosch and Qing Liu\footnote{The second author 
appreciates the hospitality of the University of M\"unster where 
this work was done.}}
\maketitle

Let $A_K$ be an abelian variety over a discrete valuation field
$K$. Let $A$ be the N\'eron model of $A_K$ over the ring of integers
$\cok$ of $K$ and $A_k$ its special fibre over the residue field
$k$ of $\cok$. Denote by ${A}^{0}$ and $A^0_k$ the corresponding
identity components. Then we have an exact sequence
$$ 0 \to {A}^{0}_k \to {A}_k \to \phi_A
\to 0 , $$
where $\phi_A$ is a finite \'etale group scheme over $k$. 
The latter is called the \emph{group of components of $A$}. 
The group of rational points 
$\phi_A(k)$ counts the number of connected components of the special
fibre 
${A}_k$ which are geometrically connected. In this paper we are 
interested in ``computing'' this group and the image 
of $A_K(K)\to \phi_A(k)$. The starting point of this 
work is an e-mail of E. Schaefer to the second author. He 
convinced us of the interest in computing $\phi_A(k)$.
\smallskip

This paper is organized as follows. Section \ref{curves} deals with 
the case where $A_K$ is the Jacobian of a curve $X_K$ over $K$. Let 
${X}$ be a regular model of $X_K$ over $\cok$. Then a modified 
intersection matrix gives an explicit subgroup of $\phi_A(k)$ and the 
quotient can be controlled by some cohomology groups. The main result 
of this section is Theorem \ref{main-jac} which determines $\phi_A(k)$ 
when $k$ is finite.

In section \ref{tors}, we put together some 
classical results and general remarks about the canonical map
$A_K(K)\to \phi_A(k)$. 

In sections \ref{Lab3} and \ref{Lab4}, we assume that $K$ is complete.
First we consider algebraic
tori $T_K$ with multiplicative reduction (so $T_K$ is not 
an abelian variety in this section). Let $T$ be the N\'eron
model of $T_K$. We show in \ref{Lab3.2} that
$\phi_T(k)$ coincides with $\phi_{T_G}(k)$, where 
$T_{G, K}$ is the biggest split subtorus of $T_K$,
and that $T(\cok)/T^0(\cok)\to \phi_T(k)$ is an isomorphism.
If $T_K$ does not admit multiplicative reduction, the same
constructions lead to subgroups of finite index; cf. \ref{Lab3.3}.
Finally, we add some results on abelian varieties 
$A_K$ with semi-stable reduction, which are more or less known.
When the toric part of 
$A_k$ is split, then $\phi_A$ is constant; cf.
\ref{Lab4.3}. In general, using data coming from 
the rigid uniformization of $A_K$, we 
are able to interpret the image of 
$A_K(K)\to \phi_A(k)$; see
\ref{Lab4.4}.
\smallskip

Throughout this paper, we fix a separable closure ${k^s}$ of $k$, and we 
denote by ${\gal}$ the absolute Galois group $\Gal({k^s}/k)$ 
of $k$.

\begin{section}{Component groups of Jacobians}\label{curves}

\newcommand{\ima}{\Img\kern 2pt\alpha}
\newcommand{\imab}{\Img\kern 2pt\ovl{\alpha}}
\newcommand{\kerb}{\Ker\kern 2pt\beta}
\newcommand{\kerbb}{\Ker\kern 2pt\ovl{\beta}}
\def\inter#1#2#3{\langle #2, #3 \rangle_{#1}}
\def\HG#1{H^1(G, #1)} 
\def\H2G#1{H^2(G, #1)} 
\def\gmi#1{\Gamma_{#1}}
\newcommand{\alb}{\ovl{\alpha}}
\newcommand{\btb}{\ovl{\beta}}

In this section, we fix a connected, proper, flat and regular curve 
$X$ over $\cok$.
Let us start with some notations. Let $\Gamma_i$, $i\in I$, be the 
irreducible components of the special fibre ${X}_k$. Denote by
${\bbz}^I$ 
the free $\bbz$-module generated by the $\Gamma_i$'s. It can be 
identified canonically with the group of Weil divisors on ${X}$ with
support 
in ${X}_k$. We denote by $d_i$ the multiplicity of $\Gamma_i$ in
${X}_k$, $e_i$ 
its geometric multiplicity (see \cite{[BLR]}, Def. 9.1.3), and let 
$r_i=[k(\Gamma_i)\cap k^{\rm s} : k]$. The integer $r_i$ is also the
number 
of irreducible components of $(\Gamma_i)_{{k^s}}$. For two divisors
$V_1, V_2$ 
on ${X}$, such that at least one of them, say $V_1$, is vertical (i.e. 
contained in ${X}_k$), we denote by $V_1\cdot V_2$ their 
{\it intersection number}
$\deg_k {\cal O}_X(V_2)|_{V_1}$. When it is necessary to refer to the
ground 
field $k$, we denote this number by $\inter{k}{V_1}{V_2}$.
\smallskip

Now let us define two homomorphisms of $\bbz$-modules which are
essential
for the computing of $\phi_A(k)$. First, 
$\alpha : {\Bbb Z}^I\to {\Bbb Z}^I$ is defined by 
$$\alpha(V)=
\sum_i r_i^{-1}e_i^{-1}\inter{k}{V}{\Gamma_i} \Gamma_i$$ 
for any $V\in \bbz^I$
(see Lemma \ref{int-num} which shows that $\alpha$ really 
takes values in $\bbz^I$). Define $\beta :  \bbz^I\to \bbz$ by 
$\beta(\Gamma_i)=r_id_ie_i$. Note that $\alpha$ can be defined 
more canonically as a map ${\Bbb Z}^I\to ({\Bbb Z}^I)'$
using a suitable (not necessarily symmetric) bilinear form. But for
our purpose, this seems not be useful.

Let ${\cal O}_K^{sh}$ denote a strict henselization of $\cok$. The
residue field of ${\cal O}_K^{sh}$ is ${k^s}$. The base change
${X}\times
\Spec {\cal O}_K^{sh}\to \Spec {\cal O}_K^{sh}$ gives rise to a
regular surface with special fibre ${X}_{{k^s}}$. Let $\ovl{I}$ be a set 
indexing the irreducible components of ${X}_{{k^s}}$. We can define
similarly 
$\ovl{\alpha} : \bbz^{\ovl{I}}\to \bbz^{\ovl{I}}$ and 
$\ovl{\beta} : \bbz^{\ovl{I}} \to\bbz$. The Galois group ${\gal}$ 
acts on $\bbz^{\ovl{I}}$ via its action on ${X}_{{k^s}}$. Moreover, it
is 
not hard to check that the action of ${\gal}$ commutes with
$\ovl{\alpha}$ 
and $\ovl{\beta}$. Note that since $\phi_A$ is \'etale over $k$,
$\phi_A({k^s})=\phi_A(k^{\rm alg})$.

\begin{thm}\label{esq-ksp} {\rm (Raynaud)} \ Let ${X}$ be a proper flat 
and regular curve over $\cok$, with geometrically irreducible generic
fibre. Assume further that either $k$ is perfect or 
${X}$ has an \'etale quasi-section. 
Let $A$ be the N\'eron model of the Jacobian of $X_K$.
Then there exists a canonical exact 
sequence of ${\gal}$-modules 
\begin{equation} \label{esq-ray} 
0 \to \imab \to \kerbb \to \phi_A({k^s}) \to 0 
\end{equation}
\end{thm}

{\PfEng}. The existence and exactness of the complex as 
abstract groups are proved in \cite{[BLR]}, Theorem 9.6.1. Let us just 
explain quickly why the map $\kerbb \to \phi_A({k^s})$ commutes with the
natural action of ${\gal}$ on both sides. To do this, let us go back 
to the construction of the map $\kerbb \to \phi_A({k^s})$ as done in 
\cite{[BLR]}, Lemma 9.5.9. 

Let $P$ be the open subfunctor of $\Pic_{{X}/\cok}$ corresponding to
line 
bundles of total degree $0$, then the N\'eron model $A$ is a 
quotient of $P$ (\cite{[BLR]}, Theorem 9.5.4). 
Since our assertion concerns only the special fibre and since the 
formation of N\'eron models commutes with \'etale base change, we can 
replace $\cok$ by a henselization and thus assume that $\cok$ is
henselian. Then ${\cal O}_K^{sh}$ is Galois over $\cok$ with
group $G$. 
Let $S=\Spec{\cal O}_K^{sh}$, $Y={X}\times S$, and let $\bbz^{\ovl{I}}$ 
be identified with $D$, group of Weil divisors of $Y$ with support in 
${X}_{{k^s}}$. Denote by $\Delta_j, j\in \ovl{I}$ the irreducible
components 
of ${X}_{{k^s}}$. Consider the homomorphism  $\rho : \Pic(Y)\to D$
defined by 
$\rho({\cal L})=\sum_{j} e_j^{-1}(\deg_k {\cal L}|_{\Delta_j})\Delta_j$. 
Then $\rho^{-1}(\kerbb)=P(S)$, and $\rho$ is a homomorphism of 
$G$-modules. Since $S$ is strictly henselian, it turns out 
that $\rho|_{P(S)}$ induces an isomorphism  $P(S)/P^{0}(S)\simeq
\kerbb/\imab$, and the canonical map $P(S)/P^{0}(S)\to 
A(S)/A^{0}(S)=\phi_A({k^s})$ is also an isomorphism. The composition 
of these two isomorphisms gives rise to the above exact sequence.
Now it is clear that all the maps we considered are compatible with
the natural actions of ${\gal}$. \qed
Now we want to take the long exact sequence of Galois cohomology of 
(\ref{esq-ray}). 
For this purpose, we need some informations on the action of ${\gal}$ on 
$\bbz^{\ovl{I}}$. Let us first state a technical lemma.

\begin{lem}\label{int-num} Let $Y$ be a proper regular scheme (of
arbitrary 
dimension) over a discrete valuation ring $\cok$. Let ${\cal O}_{K'}$ 
be a finite \'etale Galois extension of $\cok$, with residue field $k'$
at some maximal ideal of ${\cal O}_{K'}$. 
Let $\Gamma$ be an irreducible component of the special fibre $Y_k$. 
Then for any irreducible component $\Gamma'$ of $\Gamma_{k'}$, and for 
any curve $C$ (i.e. irreducible closed subscheme of dimension 1) 
contained in $Y$,
we have 
$$ \inter{k}{\Gamma}{C}=[k'\cap k(\Gamma): k]\inter{k'}{\Gamma'}{p^*C}$$
where $p$ is the canonical projection $Y_{{\cal O}_{K'}}\to Y$.
Moreover, if $Y$ has dimension $2$, then the geometric multiplicity $e$
of 
$\Gamma$ divides $\inter{k'}{\Gamma'}{p^*C}$.
\end{lem}

{\PfEng}. We consider $Y_{{\cal O}_{K'}}$ as a regular
scheme over $\cok$. Let $r=[k'\cap k(\Gamma): k]$. Then 
$p^*\Gamma=\Gamma_{k'}$ has $r$ irreducible components. Since 
$Y_{{\cal O}_{K'}}\to Y$ is Galois, and $p^*\Gamma$ is invariant by
$\Gal(K'/K)$, for any irreducible component $\Gamma''$ of $\Gamma_{k'}$
we have $\inter{k'}{\Gamma'}{p^*C}=\inter{k'}{\Gamma''}{p^*C}$. Thus
$r\inter{k'}{\Gamma'}{p^*C}=\inter{k'}{p^*\Gamma}{p^*C}=
\inter{k}{\Gamma}{C}$.
The last equality holds because ${\cal O}_{K'}$ is \'etale over $\cok$.
This proves the first part of the lemma. This part can also be proved
using
the projection formula for $p$. For the second part, we notice 
that $\Gamma'$ has the same geometric multiplicity as $\Gamma$ since $k'/k$
is
separable. Thus $e$ divides $\inter{k'}{\Gamma'}{p^*C}=
\deg_{k'} {\cal O}_{Y_{{\cal O}_{K'}}}(p^*C)|_{\Gamma'}$ according to 
\cite{[BLR]}, Cor. 9.1.8.
\qed

Before going back to groups of components, let us derive the following 
consequence.

\begin{cor}\label{appli} Let $Y$ be a regular, proper scheme over a 
discrete valuation ring $\cok$. Let $\Gamma_i$, $i\in I$, be the 
irreducible components of the special fibre $Y_k$. Denote by $d_i$ 
the multiplicity of $\Gamma_i$ in $Y_k$ and set
$r_i=[{k^s}\cap k(\Gamma_i): k]$. Then for any closed point 
$P$ of the generic fibre $Y_K$, the degree $[K(P) : K]$ is divisible by 
$\gcd\{r_id_i \ | \ i\in I\}$.
\end{cor}

{\PfEng}. Let $C:=\overline{\{ P\}}$ be the Zariski closure
of $\{ P \}$ in $Y$. Then we know that 
$$[K(P) : K]=\inter{k}{Y_k}{C}=\sum_{i\in I} d_i
\inter{k}{\Gamma_i}{C}$$
The multiplicity $r_i$ can be computed on a finite Galois extension $k'/k$ 
instead of $k^s/k$. Enlarging $k'$ if necessary, the extension 
$k'/k$ can be lifted 
to an \'etale Galois extension 
${\cal O}_{K'}/\cok$ (here `lift' means that $k'$ is the residue 
field of the localization of ${\cal O}_{K'}$ at some maximal ideal). 
According to Lemma 
\ref{int-num}, this implies that $r_i$ divides $\inter{k}{\Gamma_i}{C}$.
Thus the corollary is proved. \qed

\begin{rem} \rm This corollary confirms a prediction of 
Colliot-Th\'el\`ene and Saito (\cite{[CoS]}, Remarque 3.2 (a)). 
Actually, let $I_3=\gcd\{r_id_i \ | \ i\in I\}$, and let $I_2$
be the g.c.d of $[K(P) : K]$, where $P$ varies over the closed 
points of $Y_K$ (see \cite{[CoS]}, Th\'eor\`eme 3.1). 
Then Corollary \ref{appli} is just the divisibility relation 
$I_3\, |\, I_2$. 
We understood that in a forthcoming preprint, they will prove
that $I_1=I_2=I_3$ for $p$-adic fields. 
We think that Corollary \ref{appli} should still hold if one replaces
$r_i$ by $[k^{\rm alg}\cap k(\Gamma_i) : k]=r_ie_i$.
\end{rem}

\begin{cor} \label{div} Let ${X}$ be a connected, proper, flat and 
regular curve 
over $\cok$. Let $g$ be the genus of the generic fibre $X_K$ and let 
$d'=\gcd\{ r_id_i \ | \ i\in I\}$. Then $d' \divides 2g-2$.
\end{cor}

{\PfEng}. Let ${\cal O}_{K'}$ be a finite \'etale Galois extension of 
$\cok$ with a residue field $k'$ containing ${k^s}\cap k(\Gamma_i)$
for all $i\in I$ (see the proof of Corollary \ref{appli}). 
Denote by $p : {X}_{{\cal O}_{K'}}\to {X}$ the
projection.
Let $\omega_{{X}/\cok}$ be the relative dualizing sheaf of ${X}$. 
Consider the divisor $V:=\sum_{i\in I} {{r_id_i}\over{d'}}\Gamma_{i,0}$
on $X_{{\cal O}_{K'}}$, where $\Gamma_{i,0}$ is an irreducible component 
of $(\Gamma_i)_{k'}$. Then we have 
$$ \inter{k'}{p^*\omega_{{X}/\cok}}{V}=
\sum_{i\in I} {{d_i}\over{d'}}\inter{k}{\omega_{{X}/\cok}}{\Gamma_i}=
{{1}\over{d'}}\inter{k}{\omega_{{X}/\cok}}{{X}_k}={{2g-2}\over{d'}}
$$
This proves the corollary. 
\qed

\begin{rem}\label{smallg} \rm It is known that $d \divides g-1$
(apply the adjunction formula to ${{1}\over{d}}X_k$). It should 
be noticed that, to the contrary, $d'$ does not divide $g-1$
in general. 
\end{rem}

Now let us return to Galois action. Consider the natural injective map 
$\lambda : \bbz^I\to \bbz^{\ovl{I}}$ which sends $\Gamma$ to 
$\Gamma_{{k^s}}$. 

\begin{prop}\label{commute} Let ${X}$ be a proper flat 
and regular curve over $\cok$ with geometrically connected
generic fibre.
Then 
we have $(\bbz^{\ovl{I}})^{{\gal}}=\bbz^I$. Let $d=\gcd\{ d_i \ | \ i\in
I\}$
and $V_0:={{1}\over{d}}{X}_k$. Then we have the following commutative 
diagram of complexes:
$$
\begin{CD}
0 @>>> V_0\bbz @>>> \bbz^{\ovl{I}} @>{\ovl{\alpha}}>> \bbz^{\ovl{I}} 
@>{\ovl{\beta}}>> \bbz \\
@.      @|         @A{\lambda}AA      @A{\lambda}AA   @|  \\
0 @>>> V_0\bbz @>>> \bbz^{I} @>{{\alpha}}>> \bbz^{{I}} @>{{\beta}}>>
\bbz 
\end{CD}
$$
\end{prop}

{\PfEng}. The only non-trivial point is to check that 
$\ovl{\alpha}\circ\lambda=\lambda\circ\alpha$. Let $V\in\bbz^I$
be considered as a Weil divisor on ${X}$. Then $\alpha(V)=\sum_i 
r_i^{-1}e_i^{-1}\inter{k}{V}{\Gamma_i}\Gamma_i$. Denote by $\Gamma_{ij}$
the irreducible components of ${X}_{{k^s}}$ lying over $\Gamma_i$.
Then $\lambda(\alpha(V))=
\sum_{i, j}r_i^{-1}e_i^{-1}\inter{k}{V}{\Gamma_i}\Gamma_{ij}$. By Lemma
\ref{int-num}, $\lambda(\alpha(V))=
\sum_{i, j}e_i^{-1}\inter{{k^s}}{V_{{k^s}}}{\Gamma_{ij}}\Gamma_{ij}$.
(As in the proof of Corollary \ref{appli}, one can reduce to a finite 
Galois extension before applying Lemma \ref{int-num}). Thus
$\lambda(\alpha(V))
=\ovl{\alpha}(V_{{k^s}})$.
\qed

\begin{cor}\label{trivial}
 If $r_i=1$ (i.e. $\Gamma_i$ is geometrically irreducible) for 
all $i$, then $\phi_A$ is a constant algebraic group (or 
equivalently, $\phi_A(k)=\phi_A({k^s})$).
\end{cor}

\begin{rem} \rm It is known that for modular
curves $X_0(N)$ over ${\Bbb Q}_p$, the multiplicities $r_i$ are
equal to $1$ (at least when $N$ is square-free). Thus the component
group of the Jacobian $J_0(N)$ is constant. This fact 
was stated in \cite{[M-R]}, \S 1. 
\end{rem} 

\begin{cor} We have a canonical exact sequence of groups
\begin{equation}\label{long-seq} 
0\to \ima\to \kerb\to \phi_A(k)\to H^1({\gal}, \imab)\to H^1({\gal},
\kerbb)
\end{equation}
\end{cor}

{\PfEng}. It is clear that $(\kerbb)^{{\gal}}=\kerb$. Let us show 
that $(\imab)^{{\gal}}=\ima$. Consider the exact sequence $0\to V_0\bbz
\to \bbz^{\ovl{I}}\to \imab\to 0$, where $V_0$ is defined in the
statement of Proposition \ref{commute}, and take the long exact sequence 
of  
cohomology. It is enough to see that $H^1({\gal}, V_0\bbz)=0$. This
follows 
immediately from the facts that ${\gal}$ acts trivially on $V_0\bbz$, 
${\gal}$ is profinite and that $V_0\bbz$ has no torsion. Now we get the
corollary just by taking Galois cohomology of the exact sequence 
(\ref{esq-ray}) of Theorem \ref{esq-ksp}.
\qed

\begin{thm}\label{main-jac} Let $X$ be a proper flat and regular
curve over $\cok$ with geometrically irreducible generic fibre $X_K$.
Let 
$d=\gcd\{ d_i \ | \ i\in I\}$ and $d'=\gcd\{ r_id_i \ | \ i\in I\}$. 
Assume that $\Gal({k^s}/k)$ is procyclic (i.e. any finite extension
$k'/k$ is cyclic) and that
either $k$ is perfect, or ${X}$ has an \'etale 
quasi-section. Let $A$ be the N\'eron model of the 
Jacobian of $X_K$. Then we have an exact sequence
$$ 0 \to \kerb/\ima \to \phi_A(k) \to qd\bbz/d'\bbz \to 0$$
with $q=1$ if $d'\divides g-1$ and $q=2$ otherwise.
\end{thm}

\begin{rem}\label{elem} \rm The group $\kerb/\ima$ can be determined 
by means of elementary divisors of the matrix 
$(e_j^{-1}r_j^{-1}\inter{k}{\Gamma_i}{\Gamma_j})_{i, j\in I}$ as in
\cite{[BLR]},
Corollary 9.6.3.
\end{rem}

The remainder of the section is devoted to the proof of 
Theorem \ref{main-jac}. 

\begin{lem}\label{easy-v} Let $G'$ be a finite solvable group acting 
on a finite set $J$. Let $\bbz^J$ be endowed with the natural action 
of $G'$. Then $H^1(G', \bbz^J)=0$.
\end{lem}

{\PfEng}. First assume that $G'$ has prime order. Then
$\bbz^J$ is a direct sum of free $G'$-modules and of (free)
$\bbz$-modules
with trivial action of $G'$. Thus $H^1(G', \bbz^J)=0$. The general case 
is easily derived by induction. Note that the lemma is true for any
finite group $G'$ due to Shapira's lemma (see \cite{[Br]}, page 73).
\qed

Let $k'/k$ be a finite Galois extension containing ${k^s}\cap
k(\Gamma_i)$
for all $i\in I$. Then the components of $X_{k'}$ are geometrically
irreducible. Thus the exact sequences (\ref{esq-ray}) and
(\ref{long-seq})
can be determined over $k'$ (Corollary \ref{trivial}). For simplicity,
{\it in the rest of the proof, we denote by $G$ the group $\Gal(k'/k)$}. 
Since $G$ is cyclic, we can determine explicitly each group of 
this exact sequence.
Let us recall some notations and results of \cite{[Ser]}, VIII, \S 4. 
Fix a generator $\sigma$ of $G$. Let $m=|G|$, 
$N=\sum_{0\le j\le m-1} \sigma^j$ and $D=\sigma-1$. Recall that for any 
$G$-module $M$ we have the isomorphisms 
$$ H^1(G, M)\simeq {}_N M/DM, \quad H^2(G, M)\simeq M^G/NM$$
Moreover, if $0\to M'\to M\to M''\to 0$ is an exact sequence of
$G$-modules,
then the transition homomorphisms 
$$ \delta_1 : {}_N M''/DM'' \to M'^G/NM', \quad 
\delta_0 : M''^G \to {}_NM'/DM'$$
are given by 
\begin{equation} \label{trans}
\delta_1([x])=[Ny], \quad \delta_0([x])=[Dy]
\end{equation}
if $y\in M$ is in the preimage of $x\in M''$.

\begin{lem} \label{some-h} Recall that $V_0={{1}\over{d}}{X}_k$. The 
following properties hold :

\itboxx{i} The map $\HG{\imab}\to {{md}\over{d'}}V_0\bbz/mV_0\bbz$ 
defined by $[\ovl{\alpha}(V)]\mapsto [N(V)]$ is an isomorphism.

\itboxx{ii} Let $U\in \bbz^{\ovl{I}}$. Then $DU\in {}_N\kerbb$, and the 
map $[DU]\mapsto [\ovl{\beta}(U)]$ induces an isomorphism 
$\HG{\kerbb}\simeq d\bbz/d'\bbz$. 
\end{lem}

{\PfEng}. (i) We have an exact sequence 
$$ 0\to \HG{\imab} \to \H2G{V_0\bbz}=V_0\bbz/mV_0\bbz \to 
\H2G{\bbz^{\ovl{I}}}$$
Let $J_i$ denote the set of irreducible components of $(\Gamma_i)_{k'}$.
Then 
$$\H2G{\bbz^{\ovl{I}}}=\oplus_{i\in I} \H2G{\bbz^{J_{i}}}
=\oplus_{i\in I} \Gamma_i\bbz/mr_i^{-1}\Gamma_i\bbz$$
and the homomorphism $\H2G{V_0\bbz}\to \H2G{\bbz^{\ovl{I}}}$ sends 
$[V_0]$ to $([d_id^{-1}\Gamma_i])_i$. Then it is not hard to check 
(i) using the definition of $\delta_1$.
\smallskip

(ii) We have the exact sequence $0\to\kerbb\to \bbz^{\ovl{I}} 
\to d{\Bbb Z} \to 0$. Taking Galois cohomology we get 
$$ 0\to \Img\beta=d'{\Bbb Z} \to d\bbz\to \HG{\kerbb}\to 
\HG{\bbz^{\ovl{I}}}=0$$ 
(Lemma \ref{easy-v}). 
\qed

\noindent{\it Proof of Theorem \rm \ref{main-jac}}. Let us first 
describe the map $\psi : \HG{\imab} \to \HG{\kerbb}$ in the
exact sequence (\ref{long-seq}). One should notice that while these
groups are isomorphic, $\psi$ is not an isomorphism in general. 
Let $L : D\bbz^{\ovl{I}} \to \bbz^{\ovl{I}}$ be a section of $D :
\bbz^{\ovl{I}} \to D\bbz^{\ovl{I}}$. Let $\alb(V)\in {}_N\imab$. Since 
$\HG{\bbz^{\ovl{I}}}=0$, one has $\alb(V)\in D\bbz^{\ovl{I}}$, and thus
$\alb(V)=D(L\circ\alb(V))$. Hence using Lemma \ref{some-h} (ii), 
we see that $\psi$ is given by the formula
$$\psi([\alb(V)])=[\btb(L\circ\alb(V))]\in \HG{\kerbb}
\simeq d\bbz/d'\bbz.$$

Fix for each $i\in I$ 
an irreducible component $\gmi{i,0}$ of $(\gmi{i})_{k'}$, and put
$\gmi{i,j}:=\sigma^j(\gmi{i,0})$. Let $V_1:=\sum_i
{{r_id_i}\over{d'}}\gmi{i,0}$. 
Since $N(V_1)= {{md}\over{d'}}V_0$, Lemma \ref{some-h} (i) implies that 
$\HG{\imab}=[\alb(V_1)]d\bbz/d'\bbz$. Put 
$n:=\btb(L\circ\alb(V_1))\in\bbz$. Then $\ker\psi$ is generated by
$q[\alb(V_1)]$, where $q$ is the smallest positive integer such that 
$d'\divides qn$. Using Corollary \ref{div}, we see that to prove the 
theorem, it is enough to show that $n\equiv g-1$ mod $d'$.

Now let us construct a section of $D : \bbz^{\ovl{I}} 
\to D\bbz^{\ovl{I}}$. Since the set 
$$\{ \Gamma_{i, 0}, D\Gamma_{i,j} \  | \ i\in I, \ 0\le j\le r_i-2 \}$$
form a basis of $\bbz^{\ovl{I}}$, we have a well-defined
$\bbz$-linear map $L' : \bbz^{\ovl{I}}\to \bbz^{\ovl{I}}$ given by 
$L'(\Gamma_{i, 0})=0$, $L'(D\Gamma_{i,j})=\Gamma_{i,j}$
for any $i\in I$ and $0\le j\le r_i-2$. By construction it is clear that 
$L:=L'|_{D\bbz^{\ovl{I}}}$ is a section of $D :
\bbz^{\ovl{I}} \to D\bbz^{\ovl{I}}$. Replacing $D\Gamma_{i, j}$
by $\Gamma_{i, j+1}-\Gamma_{i, j}$, we see that 
$L'(\Gamma_{i,j})=\sum_{0\le l\le j-1}\Gamma_{i, l}$ for any
$i\in I$ and $0\le j\le r_i-1$. 
\smallskip

Let us compute the integer $n$. Applying the definitions of $\alb$ and
$L$, we get
\begin{equation}\label{compute-n}
n=\sum_{i\in I, 0\leq j\leq r_i -1} 
e_i^{-1}(V_1\cdot\Gamma_{i,j})\btb\circ L'(\Gamma_{i,j})=
\sum_{i, j} jd_i (V_1\cdot\Gamma_{i,j})=
\sum_{i} d_i(V_1\cdot U_i),
\end{equation}
where $U_i:=\sum_{0\leq j\leq r_i -1} j\gmi{i,j}$. Consider 
$W_i:=\sum_{0\leq j\leq m -1} j\gmi{i,j}$. Since $\Gamma_{i, j}=
\Gamma_{i, j'}$ if $j\equiv j'$ modulo $r_i$, we have 
(put $a=mr_i^{-1}$)
$$W_i=\sum_{0\le l\le a-1}\sum_{0\le h \le r_i-1} 
(lr_i+h)\Gamma_{i, h}
={{a(a-1)}\over{2}}r_i\gmi{i} + aU_i.$$
Since $N(V_1)={{md}\over{d'}}V_0\in V_0{\Bbb Q}$, we see that
$V_1\cdot \gmi{i}=m^{-1}(N(V_1)\cdot \gmi{i})=0$. So replacing in the
equality (\ref{compute-n}) the divisor $U_i$ by $r_im^{-1}W_i$, and 
then $V_1$ by its definition, we get 
$$ n=\sum_{i, l\in I} {{r_id_ir_ld_l}\over{md'}}\sum_{1\le j\le m-1}
j(\gmi{i,j}\cdot\gmi{l,0})$$
On the other hand, $\gmi{i,j}\cdot\gmi{l,0}=\sigma^{m-j}(\gmi{i,j})\cdot
\sigma^{m-j}(\gmi{l,0})=\gmi{i,0}\cdot\gmi{l,m-j}$. So 
$$ n=\sum_{i, l\in I} {{r_id_ir_ld_l}\over{md'}}\sum_{1\le j\le m-1} 
(m-j)
(\gmi{i,0}\cdot\gmi{l,j})$$
Adding these two equalities leads to 
$$ 2n=\sum_{i, l\in I} 
{{r_id_ir_ld_l}\over{d'}}(mr_i^{-1}\gmi{i}-\gmi{i,0})\cdot\gmi{l,0}
=\left(\sum_{l\in I} {{mr_ld_l}\over{d'}}X_{k^s}\cdot \Gamma_{l,0}
\right)-d' V_1^2
=-d' V_1^2$$
Using the adjunction formula, and with the notations of the proof of 
Corollary \ref{div}, we see that $V_1^2$ is congruent to 
$\inter{k'}{p^*\omega_{{X}/\cok}}{V_1}$ mod $2$. The latter
is equal to ${{1}\over{d'}}(2g-2)$ as calculated in the proof of 
Corollary \ref{div}. This achieves the proof of Theorem \ref{main-jac}.
\qed

\begin{rem} \rm Let $X_K$ be an elliptic curve over $K$. Let $X$ be
its minimal regular model over $\cok$. One can apply 
Theorem \ref{main-jac} and Remark \ref{elem} to compute $\phi_A(k)$. 
But one can also determine directly $\phi_A(k)$ as a subset
of $\phi_A(k^s)$ using the fact that the N\'eron model $A$ of
$X_K$ is the smooth locus of $X$. 
\end{rem}

\begin{exmp} \rm Assume that $k$ is perfect and $\chara(k)\ne 2$.
Let $a, b\in \cok$ be invertible and such that the class $\tilde{a}\in
k$
is not a square. Let $n\ge 1$ be an integer. Consider the elliptic curve 
$A_K$ given by the equation 
$$ y^2=(x^2-b\pi^{2n})(x+a) $$
where $\pi$ is a uniformizing element of $\cok$. Then the minimal
regular
model of $A_K$ over $\cok$ consists of a projective line $\Gamma_1$ over
$k$,
followed by a chain of $n-1$ projective lines over
$k(\sqrt{\tilde{a}})$, 
and ends with the conic $\Gamma_{2n}$ given by the equation 
$v^2=(u^2-\tilde{b})\tilde{a}$. Thus $\phi_A({k^s})=\bbz/2n\bbz$, 
$\phi_A(k)=\bbz/2\bbz$, and $A_k(k)/A^{0}_k(k)=\bbz/2\bbz$ or $0$
depending on $\Gamma_{2n}$ has a rational point or not. 
This shows that one cannot expect a good control of the order of 
$\phi_A({k^s})/\phi_A(k)$. 
\end{exmp}

\begin{exmp} \rm Assume $\chara(k)\ne 2$. Let $g\ge 1$, let $X_K$ be the 
hyperelliptic curve defined by an equation 
$y^2=a_0\prod_{1\le i\le g+1} (x-a_i)^2+\pi$, where $a_i\in \cok$ 
are such that their images $\tilde{a}_i\in k$ are pairwise distinct
and $\tilde{a}_0$ is not a square. Finally $\pi $ is a 
uniformizing element of $\cok$. Let $X$ be the minimal regular
model of $X_K$ over $\cok$. Then $X_k$ is integral with 
$g+1$ ordinary double points. Over $k'=k[\sqrt{\tilde{a}_0}]$, 
$X_{k'}$ splits into two components isomorphic to ${\Bbb P}^1_{k'}$
intersecting transversally at $g+1$ points. Thus using 
Theorems \ref{esq-ksp}, \ref{main-jac} and Remark \ref{elem},
we see that $\phi_A(k^s)={\Bbb Z}/(g+1){\Bbb Z}$ and $\phi_A(k)=0$.

\end{exmp}
\end{section}

\begin{section}{The homomorphism $A_K(K)\to\phi_A(k)$}\label{tors}

In this section, $A_K$ is an abelian variety over $K$. Let $A$ be 
the N\'eron model of $A_K$ over $\cok$. We would like to discuss
some relationships between $A(K)$ and $\phi_A(k)$. By the properties 
of N\'eron models, $A(\cok)=A_K(K)$. 
The specialization map gives rise to a homomorphism of groups 
$A_K(K)\to A_k(k)$. 
The second group maps canonically to $\phi_A(k)$.
In general, the map $A_k(k)\to \phi_A(k)$ is not surjective. 
The reason is that $\phi_A(k)$ counts the number of geometrically 
connected components of $A_k$, while the image of $A_k(k)$ in
$\phi_A(k)$ (which is isomorphic to $A_k(k)/A_k^0(k)$) parameterizes
the
components with rational points. Each geometrically connected component
is a torsor under ${A}^{0}_k$. But such a torsor may be non-trivial
(that is, without rational point). 

\begin{lem} Let $A_K$ be an abelian variety over $K$. 

\itboxx{i} If $K$ is henselian (e.g. complete), then 
$A_K(K)\to A_k(k)$ is surjective.

\itboxx{ii} If $k$ is finite, or if $A^{0}_k$ is an extension 
of a unipotent group by a split torus with $k$ perfect,
then $A_k(k)\to \phi_A(k)$ is surjective. 
\end{lem}

{\PfEng}. (i) Since $K$ is henselian and $A$ is smooth,
the map ${A}(\cok)\to A_k(k)$ is surjective 
(see for instance \cite{[BLR]}, Prop. 2.3.5). 

(ii) Let $k'/k$ be a finite Galois extension of $k$ such that
$A_k(k')\to \phi_A(k')$ is surjective (such an extension exists
because $\phi_A$ is finite). Then it is enough to show that 
$H^1(\Gal(k'/k), A^{0}_k(k'))=0$. 
The case $k$ finite is a theorem of Lang (\cite{[Lan]}, Theorem 2). 
The remaining case is Hilbert's 90th Theorem (see 
\cite{[Ser]}, Chap.\ X, \S 1) with induction on the dimension of 
$A^0_k$. 
\qed

\end{section}

\begin{section}{Algebraic tori}
\label{Lab3}

    In this section we consider an algebraic torus $T_{K}$ over
$K$, its N\'eron model $\tor$ over the ring of integers ${\cok}$
of $K$, and the associated component group $\phi _{\tor}$. As the
formation of N\'eron models is compatible with passing from $K$ to
its completion by \cite{[BLR]}, 10.1.3, $\phi_{\tor}$ remains
unchanged under this process, and we will assume in the following
that ${\cok}$ and $K$ are \emph{complete}. Writing 
${\cal O}_K^{sh}$ for a strict henselization of ${\cok}$ and 
$K^{sh}$ for the field of fractions
of ${\cal O}_K^{sh}$, we know then that the extension $K^{sh}/K$ is
Galois.
The attached Galois group $G$ is canonically identified with the one
of $k^{s}/k$, the residue extension of $K^{sh}/K$.

    Let us first assume that $T_{K}$ has multiplicative reduction, so
that the identity component ${\tor}_{k}^{0}$ of the special fibre 
${\tor}_{k}$ is a torus. Then $T_{K}$ splits over $K^{sh}$, and we can
view the group of characters $X$ of $T_K$ as a $G$-module. It is
well-known that in this case 
we have an isomorphism of $G$-modules
\[
\phi _{\tor} \simeq \Hom (X, {\Bbb Z}) ; 
\]
see for example \cite{[X]}, 1.1.
In particular, if $T_K$ is split over $K$, the action of $G$ on
$X$ is trivial, and $\phi_{\tor}$ is isomorphic to the constant
group ${\Bbb Z}^d$ with $d = \dim T_K$.

\begin{lem}\label{Lab3.1} Let $X_{G}$ {be the biggest
$\bbz$-free quotient of} $X$ {which is fixed by} $G$. {Then the
projection} $X \longrightarrow  X_{G}$ {gives rise to an
isomorphism}
\[
\Hom (X_{G}, {\Bbb Z}) \longrightarrow  \Hom (X, {\Bbb Z})^{G}
\]
{of groups which canonically can be identified with} 
$\phi _{\tor}(k)$, {the group of} $k$\nobd-{rational points of} 
$\phi_{\tor}$.
\end{lem}

{\PfEng}. The epimorphism $X \longrightarrow  X_{G}$ induces
injections
\[
\Hom (X_{G}, {\Bbb Z}) \hookrightarrow  \Hom (X, {\Bbb Z})^{G}
\hookrightarrow  \Hom (X, {\Bbb Z}) ,
\]
and we have to show that the left injection is, in fact, a bijection.
To do this, consider a $G$\nobd-morphism $f\colon X \longrightarrow
{\Bbb Z}$ which is fixed by $G$. Then $f$ factors through a
$G$\nobd-morphism $X/W \longrightarrow  {\Bbb Z}$ where $W \subset  X$
is the submodule generated by all elements of type $x - \sigma (x)$
with $x \in  X$ and $\sigma  \in  G$. As $X_{G}$ is obtained from
$X/W$ by dividing out its torsion part and as ${\Bbb Z}$ is
torsion-free, we see that $f$ must factor through $X_{G}$. Hence, the
map $\Hom (X_{G}, {\Bbb Z}) \hookrightarrow  \Hom (X, {\Bbb Z})^{G}$ is
bijective, as claimed. \qed
 
    Now let $T_{G,K}$ be the torus with group of characters $X_{G}$.
The projection $X \longrightarrow  X_{G}$ defines $T_{G,K}$ as the
biggest subtorus of $T_{K}$ which is split over $K$, and we can
identify the associated morphism $\Hom(X_G,{\Bbb Z})
\longrightarrow \Hom(X,{\Bbb Z})$ with the corresponding morphism of
component groups $\phi_{T_G} \longrightarrow \phi_T$. Thereby we can
conclude from \ref{Lab3.1}:

\begin{prop}\label{Lab3.2}{Let} $T_{K}$ {be a torus with
multiplicative reduction}, {and let} $T_{G,K}$ {be the
biggest subtorus which is split over $K$}. Assume that $K$ is complete.
{Then the injection} $T_{G,K} \hookrightarrow T_{K}$ {and
the associated morphism of N\'eron models} $\tor_G \longrightarrow
\tor$ induce a monomorphism of component groups $\phi _{\tor_G}
\hookrightarrow \phi _{\tor}$ and 
an isomorphism $\phi _{\tor_G}(k)
\simeq \phi _{\tor}(k)$ {between groups of}
$k$\nobd-rational points.

    Furthermore, the canonical map $T_K(K) \longrightarrow
\phi_{\tor}(k)$ is surjective, as the same is true for the split
torus $T_{G,K}$.
\end{prop}

    What can be said if, in the situation of \ref{Lab3.2}, $T_{K}$
does not have multiplicative reduction? In this case we can still
view the group of characters $X$ of $T_{K}$ as a Galois module under
the absolute Galois group of $K$. Similarly as above, we can use the
inertia group $I$ and look at the biggest subtorus $T_{I,K} \subset
T_{K}$ which splits over the maximal unramified extension $K^{sh}$ of
$K$. We get an exact sequence of tori
\[
0 \longrightarrow  T_{I,K} \longrightarrow  T_{K} \longrightarrow
\tilde T_{K} \longrightarrow  0
\]
with a torus $\tilde T_{K}$ such that $\tilde T_{K} \otimes _{K}
K^{sh}$ does not admit a subgroup of type ${\Bbb G}_{{\rm m}}$. The
N\'eron model $\tilde \tor$ of $\tilde T_{K}$ is quasi-compact by
\cite{[BLR]}, 10.2.1, and, hence, the component group $\phi _{\tilde
\tor}$ must be finite.

    We view now N\'eron models as sheaves with respect to the \'etale
(or smooth) topology on ${\cal O}_K$. Then the above exact sequence
of tori induces a sequence of N\'eron models
\[
0 \longrightarrow  \tor_{I} \longrightarrow  \tor \longrightarrow
\tilde \tor \longrightarrow  0
\]
which is exact by \cite{[BX]}, 4.2. Furthermore, using the right
exactness of the formation of component groups, see \cite{[BX]},
4.10, in conjunction with the facts that $T_{I,K}$ has multiplicative
reduction and, hence, that the component group $\phi _{\tor_{I}}$
cannot have torsion, we get an exact sequence of component groups
\[
0 \longrightarrow  \phi _{\tor_{I}} \longrightarrow  \phi _{\tor}
\longrightarrow  \phi _{\tilde \tor} \longrightarrow  0 .
\]
Restriction to $k$\nobd-rational points preserves the exactness,
\[
0 \longrightarrow  \phi _{\tor_{I}}(k) \longrightarrow  
\phi_{\tor}(k) \longrightarrow  \phi _{\tilde \tor}(k)
\longrightarrow 0 ,
\]
as $H^1(G',{\Bbb Z}^d) = \Hom(G',{\Bbb Z}^d) = 0$ for any finite
group $G'$ acting trivially on ${\Bbb Z}^d$. Now, taking into
account that $T_{I,K}$ has multiplicative reduction and that
$\phi_{\tilde \tor}(k)$ is finite, we can conclude from
\ref{Lab3.2}:

    \begin{cor}\label{Lab3.3}Let $T_{K}$ be an algebraic torus,
let $T_{G,K}$ be the biggest subtorus which is split over $K$, and
let $\tilde T_{K}$ be defined as above. Assume that $K$ is
complete. Then the canonical sequence
\[
0 \longrightarrow  \phi _{\tor_{G}}(k) \longrightarrow  
\phi_{\tor}(k) \longrightarrow  \phi _{\tilde \tor}(k)
\longrightarrow 0 ,
\]
is exact with $\phi _{\tor_{G}}(k)$ being free and $\phi _{\tilde
\tor}(k)$ finite.

	In particular, the image of $T_{G,K}(K)$ is of finite index in
$\phi _{T}(k)$, and the same is true for the image of $T_{K}(K)$.
\end{cor}

\section{Abelian varieties with semi-stable reduction}
\label{Lab4}

    Let $A_{K}$ be an abelian variety over the base field $K$, which
is assumed to be \emph{complete}. We will view $A_{K}$ as a rigid
$K$\nobd-group and use its uniformization in the sense of rigid
geometry; cf. \cite{[R]} and \cite{[BX]}, Sect.~1. So $A_{K}$ can be
expressed as a quotient $E_{K}/M_{K}$ of rigid $K$\nobd-groups with
the following properties:

\smallskip
\itboxx{i} $E_{K}$ is a semi-abelian variety sitting in a short exact
sequence
\[
0 \longrightarrow  T_{K} \longrightarrow  E_{K} \longrightarrow
B_{K} \longrightarrow  0 ,
\]
where $T_{K}$ is an algebraic torus and $B_{K}$ an abelian variety
with potentially good reduction.

\smallskip
\itboxx{ii} $M_{K}$ is a lattice in $E_{K}$ of maximal rank; i.~e., a
closed {analytic} subgroup of $E_{K}$ which, after finite
separable extension of $K$, becomes isomorphic to the constant group
${\Bbb Z}^{d}$ with $d = \dim  T_{K}$.

\smallskip
Let $\ca A$ be the N\'eron model of $A_{K}$ and $\ca A^{0}$ its
identity component. Recall that $A_{K}$ is said to have {\it semi-stable
reduction} if the special fibre $\ca A_{k}^{0}$ of $\ca A^{0}$ is
semi-abelian. Furthermore, let us talk about a \emph{split}
semi-stable reduction if the toric part of $\ca A_{k}^{0}$ is split
over $k$. The property of semi-stable reduction is reflected on the
uniformization of $A_{K}$ in the following way:

    \begin{prop}\label{Lab4.1}{The abelian variety} $A_{K}$
{has semi-stable} ({resp}. {split semi-stable})
{reduction over} $K$ {if and only if the following hold}:

\itboxx{i} {The torus} $T_{K}$ {splits over a finite
unramified extension of} $K$ ({resp}. {over} $K$).

\itboxx{ii} {The abelian variety} $B_{K}$ {has good
reduction over} $K$.

    {If the above conditions are satisfied with} $T_{K}$
{being split over} $K$, {the same is true for the lattice}
$M_{K} \subset  E_{K}$; {i}. {e}., $M_{K}$ {is then
isomorphic to the constant} $K$\nobd-{group} ${\Bbb Z}^{d}$,
{where} $d = \dim  T_{K}$.
\end{prop}

    {\PfEng}. As any abelian variety with semi-stable reduction
acquires split semi-stable reduction over a finite unramified
extension of $K$, we need only to consider the case of split
semi-stable reduction. So assume that $A_{K}$ has split semi-stable
reduction. Then we have an exact sequence
\[
0 \longrightarrow  T_{k} \longrightarrow  \ca A_{k}^{0}
\longrightarrow  B_{k} \longrightarrow  0 ,
\]
where $T_{k}$ is a split torus and $B_{k}$ an abelian variety over
$k$. Let $\cal A$ be the formal completion of $A$ along $A_k$ and
${\cal A}^0$ its identity component.
Using the infinitesimal lifting property of tori, see \cite{[SGA
3]}, exp.~IX, 3.6, and working in terms of formal N\'eron models in
the sense of \cite{[BS]}, we see that $T_{k}$ lifts to a split formal
subgroup torus ${\fltor} \subset  {\cal A}^{0}$  such
that the quotient ${\flb} = {\fla}^{0}/{\fltor}$ is a formal 
abelian scheme
lifting $B_{k}$. The theory of uniformizations, as explained for
example in \cite{[BL]}, Sect.~1, says now that the exact sequence
\[
0 \longrightarrow  {\fltor} \longrightarrow  {\fla}^{0} \longrightarrow 
{\flb}
\longrightarrow  0 ,
\]
coincides with the one obtained from
\[
0 \longrightarrow  T_{K} \longrightarrow  E_{K} \longrightarrow
B_{K} \longrightarrow  0
\]
by passing to identity components of associated formal N\'eron
models. As the group of characters of $T_{K}$ coincides with the one
of ${\fltor}$, we see that $T_{K}$ is a split torus. Furthermore, 
$\flb$ is
algebraizable with generic fibre $B_{K}$ and, thus, $B_{K}$ has good
reduction over $K$.

    Let us show that in this situation $M_{K}$ will be constant.
Indeed, writing $K^{s}$ for a separable closure of $K$, we choose
free generators of the group of characters of $T_{K}$ and look at the
associated ``valuation''
\[\begin{CD}
\nu \colon E_{K}(K^{s}) @>>> |K^{s}|^{d} @>-\log >> {\Bbb R}^{d} ,
\end{CD}\]
where $d = \dim  T_{K}$. One knows that $M_{K}$ being a lattice (of
maximal rank) in $E_{K}$ means that $M_{K}$ is of dimension zero and
that $M_{K}(K^{s})$ is mapped bijectively under $\nu $ onto a lattice
(of maximal rank) in ${\Bbb R}^{d}$.

    Now let us look at the action of the absolute Galois group 
$G_K:=\Gal(K^s/K)$ of 
$K$ on $M_{K}(K^{s})$ and show that $M_{K}$ is constant. As $K$ is
complete, the action of $G_K$ is trivial on $|K^{s}|$. Hence, it
respects the map $\nu $. Therefore $\nu $ can only be injective if
the action of $G_K$ on $M_{K}(K^{s})$ is trivial. However, then all
points of $M_{K}$ must be rational, and $M_{K}$ is constant.

    The converse, that conditions (i) and (ii) imply semi-stable
reduction for $A_{K}$, follows from \cite{[BX]}, 5.1. \qed
 
    Let us consider now an abelian variety $A_{K}$ with semi-abelian
reduction and with uniformization given by the exact sequence
\[
0 \longrightarrow  M_{K} \longrightarrow  E_{K} \longrightarrow
A_{K} \longrightarrow  0 .
\]
Then, by \ref{Lab4.1}, $M_{K}$ becomes constant over an unramified
extension of $K$, and the associated sequence of formal N\'eon models
\[
0 \longrightarrow  {\cal M} \longrightarrow  {\cal E}
\longrightarrow 
{\cal A} \longrightarrow  0
\]
is exact due to \cite{[BX]}, 4.4. As the component group $\phi _{M}$
is torsion-free, and as the formation of component groups is
right-exact, see \cite{[BX]}, 4.10, the induced sequence
\[
(*)\hskip4cm 0 \longrightarrow  \phi _{M} \longrightarrow  \phi _{E}
\longrightarrow  \phi _{A} \longrightarrow  0 \hskip4cm\strut
\]
is exact, so that $\phi _{A}$ may be identified with the quotient
$\phi_{E}/\phi _{M}$. Thus, if we view the objects of the latter
sequence as Galois modules under $G = \Gal (K^{sh}/K)$ and apply
Galois cohomology, we see:

    \begin{lem}\label{Lab4.2}As before, let $A_{K}$
be an abelian variety with semi-stable reduction. Then
the uniformization of $A_{K}$, in particular, the above
sequence \emph{($*$)}, gives rise to an exact sequence
\[
0 \longrightarrow  \phi _{M}(k) \longrightarrow  \phi _{E}(k)
\longrightarrow  \phi _{A}(k) \longrightarrow  H^{1}(G,M_{K})
\longrightarrow  \ldots 
\]
If $A_{K}$ has split semi-stable reduction, $M_{K}$
is constant and, hence, $H^{1}(G,M_{K})$ is trivial.
\end{lem}

    Thus, the quotient $\phi _{E}(k)/\phi _{M}(k)$ may be viewed as a
subgroup of the group of $k$\nobd-rational points of $\phi _{A}$, and
it coincides with $\phi _{A}(k)$ in the case of split semi-stable
reduction.

    Let $T_{K}$ be the toric and $B_{K}$ the abelian part of $E_{K}$.
Then we have an exact sequence
\[
0 \longrightarrow  T_{K} \longrightarrow  E_{K} \longrightarrow
B_{K} \longrightarrow  0
\]
of algebraic $K$\nobd-groups and, associated to it, a sequence of
N\'eron models
\[
0 \longrightarrow  T \longrightarrow  E \longrightarrow  B
\longrightarrow  0 .
\]
In terms of sheaves for the \'etale (or smooth) topology on
$\cok$, the latter is exact due to \cite{[BX]}, 4.2, as $A_{K}$
having semi-abelian reduction implies that $T_{K}$ splits over an
unramified extension of $K$; use \ref{Lab4.1} and \cite{[BX]},
5.1. Similarly as before, we get an exact sequence of component
groups
\[
0 \longrightarrow  \phi _{T} \longrightarrow  \phi _{E}
\longrightarrow  \phi _{B} \longrightarrow  0 ,
\]
where $\phi _{B}$ is trivial, since $B_{K}$ has good reduction.
Thus, the morphism $T \longrightarrow E$ induces an isomorphism
$\phi _{T} \longrightarrow^{~} \phi _{E}$ and, using the above
exact sequence ($*$), we can view $\phi _{A} = \phi_E/\phi_M$ as a
quotient $\phi _{T}/\phi _{M}$, although the morphism $M
\longrightarrow E$ might not factor through $T$.

    \begin{prop}\label{Lab4.3}Let $A_{K}$ be an abelian
variety with split semi-stable reduction; i. e.,
we assume that the identity component $A_{k}^{0}$ of the
special fibre of the N\'eron model $A$ of $A_{K}$ is
extension of an abelian variety by a split algebraic torus.
Then:

\itboxx{i} The component group $\phi _{A}$ is constant
(also valid if $K$ is not necessarily complete).

\itboxx{ii} The canonical map $A_{K}(K) \longrightarrow  \phi
_{A}(k)$ is surjective.
\end{prop}

    {\PfEng}. It follows from \ref{Lab4.1} that $M_{K}$ is constant
and that $T_{K}$ is split. Thus, the $k$\nobd-groups $\phi _{T}$ and
$\phi _{M}$ are constant, and so is their quotient $\phi _{A}$. If
$K$ is not complete, we may pass to the completion of $K$ without
changing the reduction of $A_{K}$ and its type. This establishes
assertion~(i). Furthermore, assertion~(ii) is due to the fact that
the map $T_{K}(K) \longrightarrow  \phi _{A}(k)$ is surjective, as
$T_{K}$ is a split torus. \qed

    If the semi-stable reduction of $A_{K}$ is not necessarily split,
the quotient $\phi _{T}(k)/\phi _{M}(k)$ will, in general, be a
proper subgroup of $\phi _{A}(k)$; its index is controlled by the
cohomology group $H^{1}(G,M_{K})$. To make this subgroup more
explicit, let $X$ be the group of characters of the toric part
$T_{K}$ of $E_{K}$. As is explained in \cite{[BL]}, Sect.~3 or
\cite{[BX]}, Sect.~5, we can evaluate characters of $X$ on the
lattice $M_{K}$, thereby obtaining a ``bilinear map'' $M_{K} \times
X \longrightarrow  P_{K}$ taking values in the Poincar\'e bundle
$P_{K}$ on $B_{K} \times  B^{\prime }_{K}$, the product of $B_{K}$
with its dual. In fact, if the abelian part $B_{K}$ of $E_{K}$ is
trivial, this pairing is just the evaluation of characters on the
lattice $M_{K}$. Using the canonical valuation $P_{K}(K^{sh})
\longrightarrow  \bbz$, we get a non-degenerate pairing $M_{K}
\times  X \longrightarrow  \bbz$ of Galois modules with respect to
the extension $K^{sh}/K$ and from it an injection $i\colon M_{K}
\hookrightarrow  \Hom (X,\bbz)$ into the linear dual of $X$. Now, as
a $G$\nobd-module, we can identify $M_{K}$ with the component group
$\phi _{M}$. Furthermore, $\Hom (X,\bbz)$ can be viewed as the
component group of $T_K$, and it follows from the discussion in
\cite{[BX]}, 5.2, that under this identification the inclusion map
$i$ corresponds to the canonical map $\phi _{M} \longrightarrow  \phi
_{T}$ as considered above. In particular, we have a canonical
commutative diagram of $G$\nobd-modules
\[\begin{CD}
0 @>>> \phi _{M} @>>> \phi _{T} @>>> \phi _{A} @>>> 0 \\
@. @VVV @VVV @| \\
0 @>>> M_{K} @>>> \Hom (X,\bbz) @>>> \phi _{A} @>>> 0
\end{CD}\]
with exact rows and vertical isomorphisms. Restricting to
$G$\nobd-invariants we get from \ref{Lab3.1}:

    \begin{prop}\label{Lab4.4}Let $A_{K}$ be an abelian variety
with semi-stable reduction and with uniformization $A_{K} =
E_{K}/M_{K}$. Let $X$ be the group of characters of the toric part
of $E_{K}$. Then $\Sigma = \Hom (X_{G},\bbz)/M_{K}^{G}$ is a
subgroup of $\phi _{A}(k)$, contained in the image of $A_{K}(K)
\longrightarrow \phi _{A}(k)$, such that the quotient $\phi
_{A}(k)/\Sigma$ is mapped injectively into $H^1(G,M_K)$.
Furthermore, $\Sigma$ coincides with $\phi _{A}(k)$ if $A_{K}$ has
split semi-stable reduction.
\end{prop}

    If the abelian variety does not admit semi-stable reduction, we
still have maps
\[
\phi _{T_{G}} \longrightarrow  \phi _{T_{I}} \longrightarrow  \phi
_{T} \longrightarrow  \phi _{E} \longrightarrow  \phi _{A} ,
\]
where $T_{G,K}$ stands for the maximal subtorus of $T_{K}$ which is
split over $K$ and, likewise, $T_{I,K}$ for the maximal subtorus of
$T_{K}$ which splits over $K^{sh}$. The image in $\phi _{A}$ of each
of these groups gives rise to a subgroup of $\phi _{A}$, and we
thereby get a filtration of $\phi _{A}$. Up to the term $\phi
_{T_{G}}$, this filtration was dealt with in \cite{[BX]}, Sect. 5; it
goes back to Lorenzini \cite{[L]}. Subsequent factors of the
filtration are controlled by suitable first cohomology groups or by the
component group of $B_{K}$; cf. \cite{[BX]}, 5.5. So, to make general
statements about $k$\nobd-rational points seems to be a
little bit out of reach. However, the groups $\phi _{T_{G}}$ and
$\phi _{T_{I}}$ are accessible, and this leads to the understanding
of rational components in the semi-stable reduction case.

\end{section}

\vspace{.25in}
\footnotesize
\begin{tabbing}

Mathematisches Institut \hspace{3cm} \= CNRS, Laboratoire de
Math\a'{e}matiques Pures\\
Westf\"alische Wilhelms-Universit\"at \> Universit\a'{e} de Bordeaux 1\\
Einsteinstra\ss e 62 \> 351, cours de la Lib\a'{e}ration \\
D-48149 M\"unster \> F-33405 Talence \\
bosch@math.uni-muenster.de  \> liu@math.u-bordeaux.fr \\
\end{tabbing}
\vspace{.25in}

\end{document}